\newif\ifarxiv
\arxivtrue
\ifarxiv
\documentclass{article}
\usepackage[nosumlimits]{amsmath}
\usepackage{amssymb,amsthm}
\author{C. J. Cotter, D. A. Ham, A.T.T. McRae, L. Mitchell and A. Natale}
\usepackage[width=18cm]{geometry}
\else
\documentclass[times]{qjrms4}
\runningheads{C. J. Cotter et al.}{Shallow atmosphere finite elements}
\address{<\affilnum{a}Department of Mathematics, Imperial College London\\
\affilnum{b}Department of Computing, Imperial College London\\
\affilnum{c}Grantham Institute - Climate Change and Environment, Imperial College London}
\corraddr{Department of Mathematics, Imperial College London, South Kensington Campus, SW7 2AZ, UK. E-mail: colin.cotter@imperial.ac.uk}
\author{C. J. Cotter\corrauth\affilnum{a}, D. A. Ham\affilnum{a,b}, A.T.T. McRae\affilnum{a,c}, L. Mitchell\affilnum{a,b} and A. Natale\affilnum{a}}
\fi
\usepackage{graphicx}
\def\MM#1{\boldsymbol{#1}}
\newcommand{\pp}[2]{\frac{\partial #1}{\partial #2}}

\def\MM#1{\boldsymbol{#1}}
\DeclareMathOperator{\diff}{d}

\DeclareMathOperator{\CG}{CG}
\DeclareMathOperator{\DG}{DG}
\DeclareMathOperator{\BDM}{BDM}

\DeclareMathOperator{\tri}{triangle}
\DeclareMathOperator{\interval}{interval}
\usepackage{amscd}
\usepackage{natbib}
\usepackage{fancybox}
\usepackage{tikz}
\usetikzlibrary{calc,fadings,decorations.pathreplacing}
\newcommand\pgfmathsinandcos[3]{%
  \pgfmathsetmacro#1{sin(#3)}%
  \pgfmathsetmacro#2{cos(#3)}%
}
\newcommand\LongitudePlane[3][current plane]{%
  \pgfmathsinandcos\sinEl\cosEl{#2} 
  \pgfmathsinandcos\sint\cost{#3} 
  \tikzset{#1/.estyle={cm={\cost,\sint*\sinEl,0,\cosEl,(0,0)}}}
}
\newcommand\LatitudePlane[3][current plane]{%
  \pgfmathsinandcos\sinEl\cosEl{#2} 
  \pgfmathsinandcos\sint\cost{#3} 
  \pgfmathsetmacro\yshift{\cosEl*\sint}
  \tikzset{#1/.estyle={cm={\cost,0,0,\cost*\sinEl,(0,\yshift)}}} %
}
\newcommand\DrawLongitudeCircle[2][1]{
  \LongitudePlane{\angEl}{#2}
  \tikzset{current plane/.prefix style={scale=#1}}
  \pgfmathsetmacro\angVis{atan(sin(#2)*cos(\angEl)/sin(\angEl))} %
  \draw[current plane,ultra thin,black!80] (\angVis:1) arc (\angVis:\angVis+180:1);
}
\newcommand\DrawLatitudeCircle[2][1]{
  \LatitudePlane{\angEl}{#2}
  \tikzset{current plane/.prefix style={scale=#1}}
  \pgfmathsetmacro\sinVis{sin(#2)/cos(#2)*sin(\angEl)/cos(\angEl)}
  \pgfmathsetmacro\angVis{asin(min(1,max(\sinVis,-1)))}
  \draw[current plane,ultra thin,black!80] (\angVis:1) arc (\angVis:-\angVis-180:1);
}

\tikzset{%
  >=latex,
  inner sep=0pt,%
  outer sep=2pt,%
  mark coordinate/.style={inner sep=0pt,outer sep=0pt,minimum size=3pt,
    fill=black,circle}%
}
\begin{document}

\title{On the shallow atmosphere approximation in finite element dynamical cores}

\ifarxiv
\maketitle
\fi
\begin{abstract}
  We provide an approach to implementing the shallow atmosphere
  approximation in three dimensional finite element discretisations
  for dynamical cores. The approach makes use of the fact that the
  shallow atmosphere approximation metric can be obtained by writing
  equations on a three-dimensional manifold embedded in $\mathbb{R}^4$
  with a restriction of the Euclidean metric. We show that finite
  element discretisations constructed this way are equivalent to the
  use of a modified three dimensional mesh for the construction of
  metric terms. We demonstrate our approach via a convergence test for
  a prototypical elliptic problem.
\end{abstract}
\ifarxiv
\else
\maketitle
\fi

\section{Introduction}

The shallow atmosphere approximation is a commonly used simplifying
approximation in the development of atmosphere and ocean dynamical
cores \citep{Ph1966,WhHoRoSt2005} where factors of
$1/r$ are replaced by $1/a$, where $r$ is the radial coordinate and
$a$ is some reference value (typically the Earth's mean
radius). Additionally, some of the metric terms in the momentum equation
are neglected. The approximation amounts to neglecting the increase
with $r$ of the surface area of a spherical shell with radius $r$.
\citet{Mu1989} (see also \cite{WhHoRoSt2005}) showed that the shallow
atmosphere approximation is equivalent to solving the equations of
motion in a geometry with a non-Euclidean metric.  \citet{ThWh2013}
showed that this metric can be obtained by considering the restriction of the four-dimensional Euclidean
metric to the three-manifold embedded in four-dimensional space that consists of the
surface of a sphere extruded into the fourth direction (they also
discussed the interesting aspect of non-unique geodesics in this
geometry which cause difficulties for computing Lagrangian
trajectories). We shall make use of this viewpoint in this paper. The
shallow atmosphere approximation is usually combined with the
traditional approximation, which neglects the horizontal component of
the rotation vector in the Coriolis term. This combination has historically
been used to obtain an equation set with a conserved
energy and potential vorticity; we note that an energy and PV
conserving formulation of the shallow
atmosphere approximation with non-traditional Coriolis term has
recently been discovered \citep{ToDu2014}.

In recent years there has been much interest in using finite element
methods for developing atmospheric dynamical cores
\citep{Ul2014,NaChTu2009,De+2011}, mostly because they can avoid the
parallel scalability problems associated with a latitude-longitude
grid.  In particular, the NERC/UK Met Office/STFC dynamical core
project, nicknamed ``Gung Ho'', is considering compatible/mimetic
mixed finite element methods \citep{CoSh2012}, which serve as an
extension of the C-grid staggered finite difference method, as the
basis for dynamical core development.  One option is to use a finite
element method for the horizontal discretisation and a conventional
finite difference discretisation for the vertical discretisation, with
either Lorenz or Charney-Phillips staggering. However, it is also
attractive to consider fully three-dimensional finite element
discretisations (using, for example, prismatic elements arranged in
columns) since they allow more flexibility with terrain-following
meshes.

Whilst several dynamical cores (the UK Met Office Unified Model, for
example) do not use the shallow atmosphere approximation in
operational mode, it is useful to include a minimally pervasive
shallow atmosphere approximation option since many of the standard
test cases use it. In addition, many other dynamical cores are shallow
atmosphere only, and so such an option is necessary for dynamical core
intercomparison studies.  Implementation of the shallow atmosphere
approximation in a three-dimensional finite element model is not
immediately straightforward: in contrast to finite difference methods
where the metric terms appear explicitly in the discretisation, finite
element methods can use Cartesian coordinates, with derivatives being
automatically computed \emph{via} a transformation of each element in
the mesh to a standard reference element where all differentiation and
integration is performed with the result that the metric terms do not
appear explicitly in the discretisation.  This leads to the question
of how to apply the shallow atmosphere approximation within a
three-dimensional finite element approximation. The solution lies in
the use of the transformation from four-dimensional space discussed in
\citet{ThWh2013}, referred to above. By first mapping the equations
onto the embedded 3-manifold in four dimensional space, and
\emph{then} transforming back to the reference element, we obtain
equations which are equivalent to the shallow atmosphere
approximation.

The rest of this paper is organised as follows. In Section
\ref{sec:pullback} we review the calculus tools associated with the
shallow atmosphere approximation. In Section \ref{sec:fem}, we explain
our approach to implementing the shallow atmosphere approximation in
finite element models. We then illustrate this approach with a
numerical convergence test in Section \ref{sec:num}, and finally
provide a summary and outlook in Section \ref{sec:summary}.

\section{Shallow atmosphere approximation through pullback}
\label{sec:pullback}
In this section, we review the necessary calculus tools to form
shallow atmosphere approximations of equations. We are interested
in solving equations in the spherical annulus domain
\begin{equation}
\begin{split}
M & = \left\{\MM{x}=(x_1,x_2,x_3)\in \mathbb{R}^3:\right. \\
 & \qquad \qquad \left. a^2\leq x^2_1+x^2_2+x^2_3\leq
 (a+H)^2\right\} \subset \mathbb{R}^3,
\end{split}
\end{equation}
where $a$ denotes the radius of the interior spherical surface and $H$
denotes the thickness of the spherical annulus. As discussed in
\cite{ThWh2013}, the shallow atmosphere approximation can be obtained
by writing equations in the domain
\begin{equation}
\begin{split}
\tilde{M} &= \left\{\tilde{\MM{x}}=(\tilde{x}_1,\tilde{x}_2,
\tilde{x}_3,\tilde{x}_4)\in \mathbb{R}^4:\right. \\
& \left. \qquad \qquad a^2= \tilde{x}^2_1+\tilde{x}^2_2+\tilde{x}^2_3,\,
0\leq \tilde{x}_4 \leq H\right\} \subset \mathbb{R}^4,\\
\end{split}
\end{equation}
and transforming to the physical domain $M$ using the
smooth invertible map
$\Phi:\tilde{M}\to M$ defined by
\begin{equation}
\label{eq:Phi}
\tilde{\MM{x}}=(\tilde{x}_1,\tilde{x}_2,\tilde{x}_3,\tilde{x}_4)
 \mapsto \MM{x} = 
\Phi(\tilde{\MM{x}}) =  \left(1+\frac{\tilde{x}_4}{a}\right)
(\tilde{x}_1,\tilde{x}_2,\tilde{x}_3).
\end{equation}
This transformation contains two aspects. Firstly, all derivatives
must be modified using the chain rule. Secondly, the various physical
quantities must be transformed as follows:
\begin{enumerate}
\item Scalar quantities (such as potential temperature) must be transformed
according to 
\begin{equation}
\tilde{\theta} \mapsto \theta = \tilde{\theta}\circ\Phi^{-1}.
\end{equation}
\item Vector fields (such
as the velocity field) must be transformed
according to 
\begin{equation}
\tilde{\MM{u}} \mapsto \MM{u}
= \left(J\circ \Phi^{-1}\right)^{-T}
\tilde{\MM{u}}\circ \Phi^{-1}, \quad J := \tilde{\nabla}\Phi,
\end{equation}
where the superscript $^{-T}$ indicates the transpose of the inverse.
Since we have a mapping from a 3-manifold $\tilde{M}$ embedded in
$\mathbb{R}^4$ to $\mathbb{R}^3$, the inverse is defined on the
tangent space to $\tilde{M}$. If the $\mathbb{R}^4$ coordinate system
is used to represent tangent vectors, then the Moore-Penrose inverse
of $J$ should be used.
\item Fluxes (\emph{i.e.}, vector quantities that must be integrated 
over surfaces to obtain flow rates, such as the mass flux) must
be transformed according to 
\begin{equation}
\tilde{\MM{F}} \mapsto \MM{F} = \frac{1}{\det J}
\left(J\circ \Phi^{-1}\right)\tilde{\MM{F}}\circ\Phi^{-1},
\end{equation}
where $\det J$ is the pseudodeterminant of $J$, namely the product
of the non-zero singular values of $J$.
\item Densities (\emph{i.e.}, scalar quantities that must be integrated
over volumes to obtain total quantities, such as kinetic energy density)
must be transformed according to
\begin{equation}
\tilde{\rho} \mapsto \rho = \frac{1}{\det J}\tilde{\rho}\circ\Phi^{-1}.
\end{equation}
\end{enumerate}
The
resulting evolution equations for the quantities in $M$ will then
contain metric terms that encode the shallow atmosphere approximation.

For example, we could take the velocity equation and write it on
$\tilde{M}$ as
\begin{align}
\label{eq:shallow atmosphere Mtilde}
\pp{\tilde{\MM{u}}}{t} + (\tilde{\MM{u}}\cdot\nabla)\tilde{\MM{u}}
+ 2\tilde{\MM{\Omega}}{\times}\tilde{\MM{u}}
 &= -\tilde{\theta}\nabla p + \nabla\tilde{\phi}
+ \MM{l}\lambda, \\
\MM{l}\cdot\tilde{\MM{u}} & = 0,
\end{align}
where $\tilde{\MM{u}}$ is the velocity field, $\tilde{\theta}$ is the
potential temperature, $\tilde{p}$ is the pressure, and $\tilde{\phi}$ is
the geopotential, all defined on $\tilde{M}$. Further, $\MM{l}$ is the
unit vector that is normal to all tangent vectors on $\tilde{M}$,
$\lambda$ is a Lagrange multiplier field that enforces the condition
that $\tilde{\MM{u}}$ stays tangential to $\tilde{M}$, and the
gradient operator $\nabla$ is restricted to the tangent space on
$\tilde{M}$.

 The cross-product of
two vector fields $\tilde{\MM{u}}$ and $\tilde{\MM{v}}$ on the
3-manifold $\tilde{M}$ is defined in the usual way as
\begin{equation}
\tilde{\MM{u}}\times\tilde{\MM{v}}
= \|\tilde{\MM{u}}\|\|\tilde{\MM{v}}\|\sin\theta\MM{n},
\end{equation}
where $\theta$ is the angle between $\tilde{\MM{u}}$ and
$\tilde{\MM{v}}$, and $\MM{n}$ is the unit vector perpendicular to
both $\tilde{\MM{u}}$ and $\tilde{\MM{v}}$ but still in the tangent
space to $\tilde{M}$, with the sign determined according to the
right-hand rule as usual. This can be computed by finding an
orthonormal basis for the tangent space at each point, expanding
$\tilde{u}$ and $\tilde{v}$ in that basis and using the standard
formula for the corresponding components of $\tilde{w}$.  We could
also make the traditional approximation, so that
\begin{equation}
\tilde{\MM{\Omega}} = \frac{\Omega x_3}{r^2}\left(x_1,x_2,x_3,0\right),
\quad r^2=x_1^2+x_2^2+x_3^2,
\end{equation}
where $\Omega$ is a scalar constant. In which case we obtain
\begin{equation}
\tilde{\MM{\Omega}}\times\tilde{\MM{u}}
= \frac{\Omega x_3}{r^2}\left(
x_2u_3 - x_3u_2,-x_1u_3 + x_3u_1,
x_1u_2 - x_2u_1,
u_4
\right).
\end{equation}
To extend these ideas to domains with varying topography, we keep the
same transformation map $\Phi$ from Equation \eqref{eq:Phi}, but
modify the corresponding domains $M$ and $\tilde{M}$.

\section{Finite element methods for the shallow atmosphere
  approximation}
\label{sec:fem}

Dynamical cores using finite difference discretisations are usually
developed by starting with the equations written in spherical polar
coordinates. This can be thought of as writing the equations (such as
\eqref{eq:shallow atmosphere Mtilde}) on $M$, applying the chain rule
to the coordinate change from $M$ into spherical coordinates, and
applying the appropriate transformations to all of the physical
quantities. This transformation leads to metric terms appearing in the
equations.  Conversely, the shallow atmosphere approximation is
obtained by writing the equations on $\tilde{M}$ and transforming from
there instead of $M$, leading to relevant alterations of the
metric terms.

In contrast, when using a finite element discretisation it is often
more natural to keep everything in Cartesian coordinates on $M$,
particularly since this then avoids the awkward problem of how to
discretise the spherical coordinate metric terms whilst maintaining
conservation, stability, \emph{etc}.  Hence, we propose to obtain
finite element discretisations of the shallow atmosphere approximation
by solving Equation \eqref{eq:shallow atmosphere Mtilde} in
$\tilde{M}$ (together with the equations of motion for the other
quantities) before finally transforming back to spherical coordinates
to produce the results. If integral quantities are required for
postprocessing, such as circulation loop integrals of velocities, flux
integrals through surfaces, or volume integrals, then the metric terms
must be included otherwise the correct conservation properties will
not be observed. The easiest way to do this is to compute these
quantities directly on $\tilde{M}$.

To begin, we define an approximation to the curved manifold
$\tilde{M}$ by selecting a discrete set of points in $\tilde{M}$ and
then using these points as vertices in a mesh.  We denote this
mesh by $\tilde{M}^\delta$.  This choice also defines a mesh
approximating $M$ in $\mathbb{R}^3$ which we denote $M^\delta :=
\Phi(\tilde{M}^\delta)$.

In contrast to finite different methods, rather than discretising the
operators in the PDE, we choose discrete function spaces defined on
$\tilde{M}^\delta$ in which we
seek a solution weakly (\emph{i.e.} in integral form).  The core of
the finite element method therefore boils down to integrating known
functions over the domain $\tilde{M}^\delta$, typically via numerical
quadrature.  These integrals may be rewritten as sums of integrals
over the individual elements making up $\tilde{M}^\delta$.  The
integration is then performed by transforming from each individual
\emph{physical} element to a \emph{reference} element $\hat{e}$.  This
merely requires that we have a transformation $g_e$ from the reference
element to each physical element $e$: gradients are then computed by
the change-of-variables formula requiring the inverse of the Jacobian
$J_e$ of $g_e$.

This extends naturally to problems solved on an n-dimensional manifold
$\mathcal{M}$ embedded in $\mathbb{R}^m, \: m > n$, as described in
\citet{RoHaCoMc2013}.  We just need to be careful since $J_e$ is no
longer square and its inverse is therefore not well-defined.  It may,
however, be inverted under the assumption that the solution lies in
the tangent space of $\mathcal{M}$ by using the Moore-Penrose inverse.

\begin{figure}[htbp]
  \centering
  \begin{tikzpicture}[scale=0.7,every node/.style={minimum size=1cm}]
  
    \def\R{4} 
  
    \def\angEl{25} 
    \def\angAz{-100} 
  
  
    \LongitudePlane[pzplane]{\angEl}{-50}
    \LongitudePlane[qzplane]{\angEl}{-60}
    \LongitudePlane[rzplane]{\angEl}{-56} \fill[ball color=white!90]
    (0,0) circle (\R); 

    \path[pzplane] (30:\R) coordinate (P1); \path[pzplane] (30:\R+3)
    coordinate (P1d); \path[rzplane] (37:\R) coordinate (P2);
    \path[rzplane] (37:\R+3) coordinate (P2d); \path[qzplane] (27:\R)
    coordinate (P3); \path[qzplane] (27:\R+3) coordinate (P3d);

    \foreach \alpha in {0, 20, ..., 160}
    \DrawLongitudeCircle[\R]{\alpha}; \foreach \beta in {-70, -50,
      ..., 70} \DrawLatitudeCircle[\R]{\beta};
  
    \draw[ultra thick, black, line cap=round, join=round] (P1) --
    (P2) -- (P3) -- (P1) -- (P1d) -- (P2d) -- (P3d) -- (P1d) -- (P1)
    (P2) -- (P2d) (P3d) -- (P3);

    \path (6, 0) coordinate (R1); \path ($(R1) + (-70:1)$) coordinate
    (R2); \path ($(R2) + (30:1.5)$) coordinate (R3); \path ($(R1) +
    (90:1)$) coordinate (R1d); \path ($(R2) + (90:1)$) coordinate
    (R2d); \path ($(R3) + (90:1)$) coordinate (R3d);

    \draw[ultra thick, black, line cap=round, join=round] (R1) --
    (R2) -- (R3) -- (R1) -- (R1d) -- (R2d) -- (R3d) -- (R1d) -- (R1)
    (R2) -- (R2d) (R3d) -- (R3);

    \path ($(R1)!0.5!(R2)$) coordinate (R12); \path
    ($(R1d)!0.5!(R2d)$) coordinate (R12d);

    \path ($(P3d)!0.5!(P1d)$) coordinate (P13d); \draw[ultra thick,
    ->] ($(R12)!0.5!(R12d) + (-0.2, -0.2)$ ) to [bend left=20]
    node[below] {$g_e$} ($(P13d) + (0.1, -0.1)$);
  \end{tikzpicture}
  \caption{Mapping from 3-dimensional reference prismatic element to a
    physical element in the spherical annulus.  Note how the physical
    element has an increasing area when moving radially outwards, such
    that $g_e$ is non-affine.}
\label{fig:mapping-3d}
\end{figure}
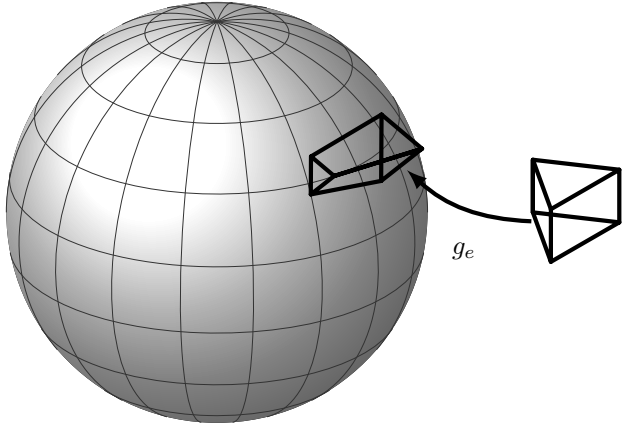

For large-scale geophysical flows it is usually desirable to use
finite element meshes arranged in vertical columns, resulting in
prismatic elements as shown in figure \ref{fig:mapping-3d}. In
contrast to tetrahedral elements, if prismatic elements are used,
meshes of the spherical annulus domain $M$ will result in the
transformation $g_e$ from the reference element to the physical
element being non-affine (affine meaning the combination of a
translation and a linear map). This is because the triangle at the top
of each prismatic element in the mesh is larger than than the triangle
at the bottom. The result is that $J_e$ is not constant within each
element and must be recalculated (and if necessary, inverted) at each
quadrature point.  In the absence of varying topography, and in the
case of elements with straight sides in $\mathbb{R}^4$, a useful
side-effect of the shallow atmosphere approximation is that this area
increase does not occur in $\tilde{M}$.  Hence $g_e$ \emph{is} affine
and $J_e$ need only be computed (and inverted) once per element,
reducing the required number of floating point operations.

\subsection{Solving in 4-d}
\label{sec:solving}
\begin{figure}[htbp]
  \centering
  \begin{tikzpicture}
    \coordinate (R1) at (0, 0);
    \coordinate (R2) at (1, 0);
    \coordinate (R3) at (0, 1);
    
    \coordinate (M1) at (2, -0.2);
    \coordinate (M2) at (3, 0.3);
    \coordinate (M3) at (2.5, 1.2);

    \coordinate (Mt1) at (0.7, 2);
    \coordinate (Mt2) at (1.9, 2.3);
    \coordinate (Mt3) at (1.3, 3.7);

    \draw[very thick, join=round, line cap=round] (R1) -- (R2) -- (R3) -- cycle;
    \node at (barycentric cs:R1=1,R2=1,R3=1) {$\hat{e}$};

    \draw[very thick, join=round, line cap=round] (M1) -- (M2) -- (M3) -- cycle;
    \node at (barycentric cs:M1=1,M2=1,M3=1) {$e$};

    \draw[very thick, join=round, line cap=round] (Mt1) -- (Mt2) -- (Mt3) -- cycle;
    \node at (barycentric cs:Mt1=1,Mt2=1,Mt3=1) {$\tilde{e}$};

    \draw [ultra thick, ->] ($(R1)!0.5!(R3) + (-0.1, 0)$) to[bend left=80] node[left, inner sep=5] {$\tilde{g}_e$} ($(Mt1)!0.5!(Mt3) + (-0.1, 0)$);

    \draw [ultra thick, ->] ($(R1)!0.5!(R2) + (0, -0.1)$) to[bend right=80] node[below, inner sep=5] {$g_e$} ($(M1)!0.5!(M2) + (0, -0.1)$);

    \draw [ultra thick, ->] ($(Mt2)!0.5!(Mt3) + (0.1, 0)$) to[bend left=80] node[right, inner sep=5] {$\chi_e := g_e \circ \tilde{g}_e^{-1}$} ($(M2)!0.5!(M3) + (0.1, 0)$);

  \end{tikzpicture}
  \caption{Transformations between the reference element ($\hat{e}$) and physical
    elements in $\mathbb{R}^4$ ($\tilde{e}$) and $\mathbb{R}^3$ ($e$)}
  \label{fig:transformations}
\end{figure}
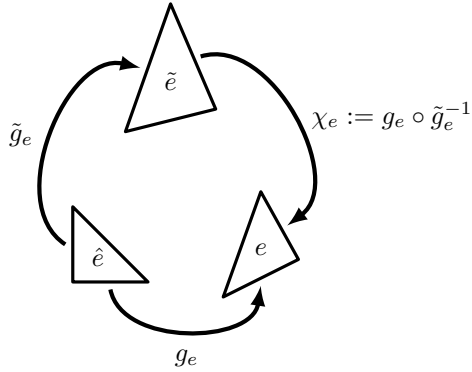

Many finite element software libraries do not contain the capability
to solve equations on embedded manifolds in higher-dimensional
spaces. Further, the Fenics 1.2 implementation of \citet{RoHaCoMc2013}
did not consider the case of a 3-manifold embedded in $\mathbb{R}^4$.
One alternative would be to solve the equations in $\mathbb{R}^3$ but
include the metric terms obtained from the transformation from
$\tilde{M}$. However, this would be a very pervasive change, since it
changes the equations at the element integral level, and therefore
limits the possibilities for using shallow atmosphere approximation
tests to build up confidence in the deep atmosphere version of a
dynamical core.  Fortunately, we can circumvent this problem by
composing two transformations: the first ($\tilde{g}_e$) from the
reference element $\hat{e} \subset \mathbb{R}^3$ into
$\tilde{M}^\delta \subset \mathbb{R}^4$; the second ($g_e$) from
$\hat{e}$ into $\mathbb{R}^3$.  The relationship between these
transformations is sketched in figure \ref{fig:transformations}.

$\tilde{g}_e : \hat{e} \to \tilde{M}^\delta$ is just the coordinate
transformation from $\hat{e}$ into the element $\tilde{e}$ in
$\tilde{M}^\delta$, expressible by expanding in a nodal basis with
coefficients given by the action of $\tilde{g}_e$ on the chosen node
points in $\hat{e}$.  $g_e : \hat{e} \to \mathbb{R}^3$ is defined as
follows:
\begin{enumerate}
\item Define the projector $P : \mathbb{R}^4 \to \mathbb{R}^4$ by 
\[
P(\tilde{\MM{x}}) := P(\tilde{x}_1,\tilde{x}_2,\tilde{x}_3,\tilde{x}_4) = (\tilde{x}_1,\tilde{x}_2,\tilde{x}_3,0),
\]
\emph{i.e.} the map leaves the first three components the same
but maps the fourth component to zero.
\item Compute the average $\tilde{\MM{x}}_e$ of the vertices of 
  the element $\tilde{e}$.
\item Compute the unit vector $\MM{k}_e \in \mathbb{R}^3$ which is 
  normal to the sphere of radius $a$ at $\Phi(P(\tilde{\MM{x}}_e))$.
\item For each node point $\MM{\xi}$ in
$\hat{e}$, define the nodal value of $g_e$ by
\begin{equation}
g_e(\MM{\xi}) := \Phi(P(\tilde{g}_e(\MM{\xi})))
+ \left(\tilde{\MM{i}}_4\cdot\tilde{g}_e(\MM{\xi})\right)\MM{k}_e
\end{equation}
where $\tilde{\MM{i}}_4 := (0, 0, 0, 1)$.
Finally, use these node points as basis coefficients to define
$g_e$ throughout the element.
\end{enumerate}
Note that while $\tilde{g}_e$ maps from the reference element into
$\tilde{M}^\delta$, $g_e$ \emph{does not} map into $M^\delta$.

We can now define the transformation $\chi:\tilde{M}^\delta \to
\mathbb{R}^3$ elementwise by
\begin{equation}
\chi_e = g_e\circ\tilde{g}_e^{-1},
\end{equation}
where $\chi_e$ is the restriction of $\chi$ to the element
$\tilde{e}$.  Note that in general $\chi_e$ is discontinuous between
elements.

Having defined $g_e$, we can see that transforming from the
reference element $\hat{e}$ into element $\tilde{e}=\tilde{g}_e(\hat{e})$ on
$\tilde{M}^\delta$ is equivalent to first transforming to element
$e=g_e(\hat{e})$, and then transforming from element $e$ to element
$\tilde{e}$ using $\chi_e^{-1}$. Equivalently, recall that the
transformations enter our equations through the Jacobian of the
mapping between domains, we can therefore compute directly in
$M^\delta$ by arranging that the Jacobians are computed correctly.  We
achieve this by using the (discontinuous, piecewise polynomial)
coordinate field defined by:
\begin{equation}
\MM{x}' := \chi \circ \Phi^{-1}(\MM{x})
\end{equation}
where $\MM{x}$ is the original coordinate field on $M^\delta$.  At
this point we no longer need to compute (or transform to) on
$\tilde{M}^\delta$ at all, instead we just compute $J_e$ on each
element $e$ in $M^\delta$ using $\MM{x}'$ rather than $\MM{x}$.  This
can be implemented in finite element software provided that it
supports discontinuous coordinate fields (this is the case in the
Firedrake software library, for example).

\section{Numerical example}
\label{sec:num}
In this section, we verify that this approach results in correctly
implementing the shallow atmosphere approximation, applied to the
prototype linear elliptic system of equations
\begin{equation}
\begin{split}
\MM{u} + 2\MM{\Omega}\times\MM{u} &= -\nabla p + \MM{F}, \\
\nabla \cdot \MM{u} - p = g, \\
\end{split}
\end{equation}
where $\MM{u}$ and $p$ are the vector and scalar unknowns
respectively, $\MM{\Omega}$ is a rotation
vector field under the traditional approximation, $\MM{e}_r$ is the
unit vector in the radial direction, $\phi$ is the latitude,
$\MM{F}$ is a prescribed vector valued forcing field and $g$ is a prescribed
scalar source. The equations are solved in a spherical annulus domain,
with inner radius $a=1$ and outer radius 2 (this leads to a very big difference
with and without shallow atmosphere approximation to demonstrate that
the method is working).
respectively.  We use these equations since we can easily construct
exact solutions,
they are linear and known to be well-posed which means that we can make
reliable convergence analyses, and we avoid technical details of
precise discretisations of nonlinear advection terms \emph{etc.} They
are also similar to equations that arise in the linear solver step in
a semi-implicit formulation of the three dimensional compressible
Euler equations. If we obtain the correct order of convergence of
numerical solutions of these equations then we will demonstrate that
the approach works correctly.

In this test case we take
\begin{eqnarray}
\MM{\Omega} &= 
\begin{pmatrix}
0 \\
0 \\
0 \\
\frac{1}{2}x_3, \\
\end{pmatrix}, \\
g & = x_1x_2x_3(x_4^2-1)(x_4^2-4), \\
\MM{F} &= 
x_3\begin{pmatrix}
\left(x_2^2
-x_3^2\right)x_1(x_4^2-1)(x_4^2-4) \\
\left(x_3^2
-x_1^2\right)x_2(x_4^2-1)(x_4^2-4) \\
\left(x_1^2
- x_2^2\right)x_3(x_4^2-1)(x_4^2-4) \\
0
\end{pmatrix},
\end{eqnarray}
where we have defined the functions on $\tilde{M}$ in terms of
coordinates in $\mathbb{R}^4$. In this case, the equations have a
unique solution given by
\begin{eqnarray}
p &= x_1x_2x_3(x_4^2-1)(x_4^2-4), \\
\MM{u} &= 
\begin{pmatrix}
x_2x_3(1-x_1^2)(x_4^2-1)(x_4^2-4) \\
x_1x_3(1-x_2^2)(x_4^2-1)(x_4^2-4) \\
x_1x_2(1-x_3^2)(x_4^2-1)(x_4^2-4) \\
2x_1x_2x_3x_4(2x_4^2-5)
\end{pmatrix}.
\end{eqnarray}

The finite element discretisation is obtained first transforming the
equations to weak form.  This is done by multiplying both equations by
test functions $\MM{w}$ and $\phi$, and integrating over the domain,
and we get
\begin{align} \nonumber
\int_{M} \MM{w}\cdot \MM{u} + 2\MM{w}\cdot\MM{\Omega}\times
\MM{u}  & \\
\qquad \qquad  - \, \nabla\cdot{\MM{w}}p\diff x &= 
\int_{M} \MM{w}\cdot\MM{F} \diff x, \\
\int_M \phi\left(\nabla\cdot\MM{u} - p\right)\diff x
&= \int_M \phi g\diff x.
\end{align}
The finite element approximation is obtained by replacing $M$ by
$M^\delta$, restricting $\MM{w}$ and $\MM{u}$ to a chosen vector
finite element space $V_1$, and restricting $\phi$ and $p$ to a scalar
finite element space $V_2$, and we get
\begin{align}
\nonumber
\int_{M^\delta} \MM{w}^\delta\cdot \MM{u}^\delta + 2\MM{w}^\delta\cdot\MM{\Omega}^\delta\times
\MM{u}^\delta & \\
\qquad \qquad -\, \nabla\cdot{\MM{w}}^\delta p^\delta\diff x &= 
\int_{M^\delta} \MM{w}^\delta\cdot\MM{F} \diff x, \quad
\forall \MM{w}^\delta \in V_1, \\
\int_{M^\delta} \phi^\delta\left(\nabla\cdot\MM{u}^\delta - p^\delta\right)\diff x
&= \int_{M^\delta} \phi^\delta g\diff x, \quad
\forall \phi^\delta \in V_2.
\end{align}
To obtain the shallow atmosphere approximation, all the factors of $J_e$
in each element integral are computed using the transformation $\hat{g}_e$
defined above, \emph{i.e.}, the discontinuous coordinate field $\MM{x}'$ 
is used to compute metric terms instead of $\MM{x}$.

\begin{figure}
\centerline{
\includegraphics[width=7cm]{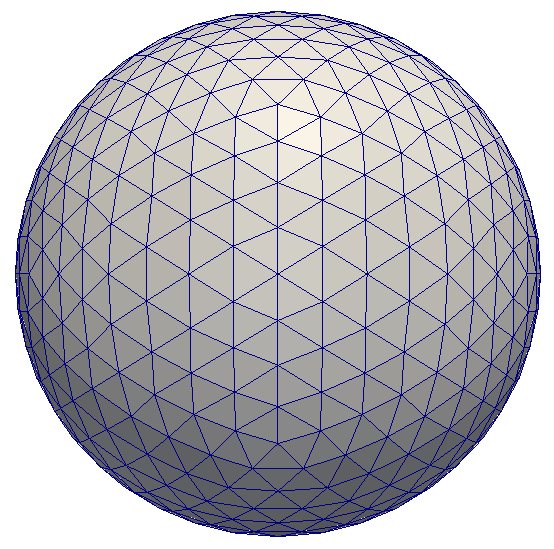}
}
\centerline{
\includegraphics[width=7cm]{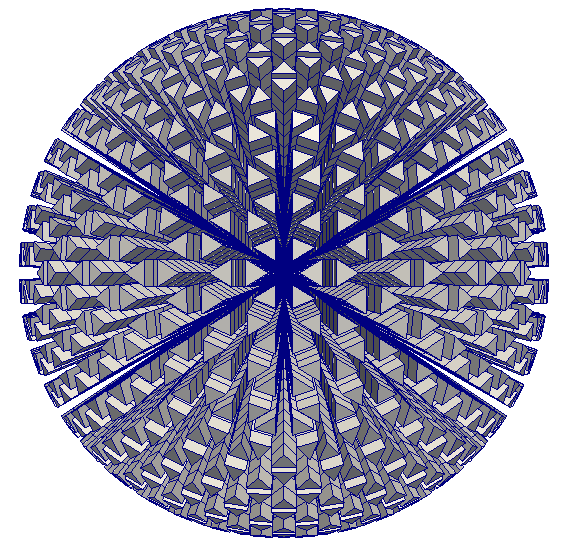}
}
\caption{\label{fig:hedgehog}A visualisation of the ``hedgehog''
  domain that is obtained by applying the mapping $\chi$ to a mesh of
  the spherical annulus constructed from columns of triangular prism
  elements. {\bfseries Top}: the original domain
  $M^\delta$. {\bfseries Bottom}: the transformed domain
  $\chi(M^\delta)$.}
\end{figure}

In this particular experiment, we use an extruded mesh for $M$, made
of triangular prism elements arranged in columns. Under the
transformation $\MM{x}$ to $\MM{x}'$ defined above, this leads to a
mesh in which each column stays the same width from bottom to top, and
hence has gaps between each column. This domain, which we have informally named  the
``hedgehog mesh'', is illustrated in Figure \ref{fig:hedgehog}. We reiterate that the discontinuous coordinate
field does not imply any loss of continuity in the underlying finite
element spaces, since the topology of the original extruded mesh is
used (and quantities are be mapped back to $M$ as a postprocessing
step).

Since our weak form of the equation only involves the divergence of
test and trial functions, we can consider using H(div) elements for
$V_1$, and discontinuous elements for $V_2$, chosen such that
\[
\MM{u}^\delta \in V_1 \implies \nabla\cdot\MM{u}^\delta \in V_2,
\]
which leads to a stable discretisation. In our test, we use the
natural extension of Brezzi-Douglas-Marini (BDM) finite element spaces
to triangular prisms, which we can express as
\[
V_1 = \BDM_k(\mbox{prism}) = 
\begin{pmatrix}
\BDM_k(\tri)\otimes \DG_{k-1}(\interval) \\
\DG_{k-1}(\tri) \otimes \CG_k(\interval) \\
\end{pmatrix},
\]
where $\DG_k$ indicates a discontinuous finite element space of degree
$k$, $\CG_k$ indicates a continuous finite element space of degree
$k$, and $\otimes$ indicates a tensor product. Here the finite element
space is described as a vector with the horizontal part of the vector
fields above and the vertical part of the vector fields below.  The
finite element spaces are defined on the reference element $\hat{e}$
and transformed \emph{via} the contravariant Piola transformation to
obtain vector-valued functions in the physical elements in
$\tilde{M}^\delta$ which are tangential to $\tilde{M}^\delta$, and
have continuous normal components across element edges, by
construction.  This means that it is not necessary to include Lagrange
multipliers to enforce the tangency constraint. Under the mapping
$\chi:\tilde{M}^\delta \to \mathbb{R}^3$, we apply a further Piola
transformation, which is equivalent to applying the Piola
transformation from the reference element into the hedgehog mesh
directly; this means that normal components of $\MM{u}^\delta$ are the
same on either side of the jump between two neighbouring columns in
$\chi(\tilde{M}^\delta)$. 

\begin{figure}
\centerline{
\includegraphics[width=10cm]{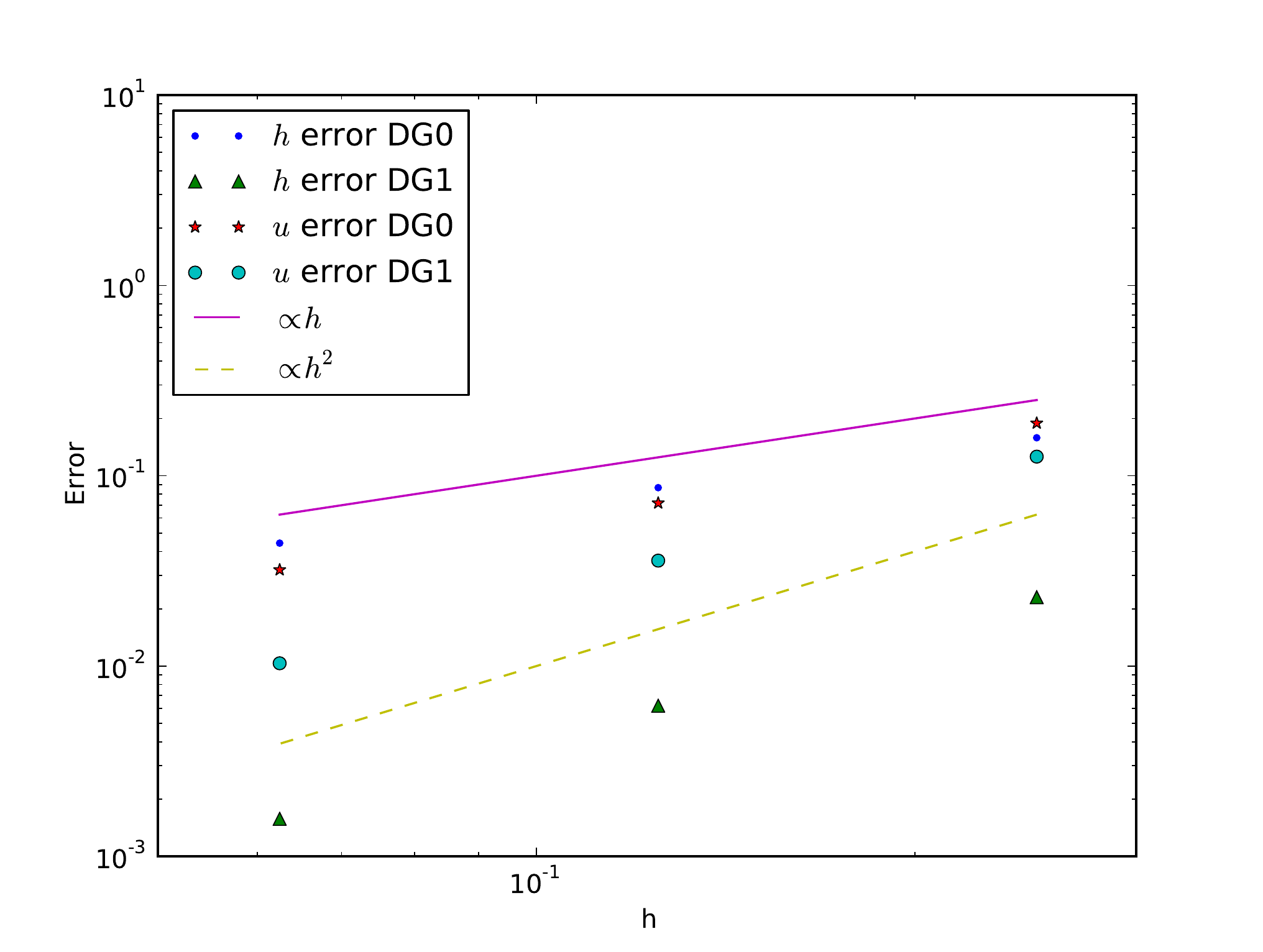}
}
\caption{\label{fig:errs}Plot showing errors in $h$ and $\MM{u}$ for
  $k=1$ and $k=2$.  For $k=1$ we obtain first order convergence for
  both $h$ and $\MM{u}$. For $k=2$ we obtain secord order convergence
  for $h$ but sub-optimal convergence (somewhere between first and
  second order) for $\MM{u}$; this is because we have used a piecewise
  linear approximation to the annular domain.}
\end{figure}

The corresponding discontinuous finite element space is
\[
V_2 = \DG_k(\mbox{prism}) =
\DG_{k-1}(\tri)\otimes \DG_{k-1}(\interval).
\]
This pair of finite element spaces satisfies the Brezzi stability
conditions with respect to our equations, and hence we expect
convergence of numerical solutions at the optimal rate, provided that
the spherical annulus domain is approximated at the correct order.
\cite{holst2012geometric} showed that this requires the correct order
of convergence for not only the distance from the approximate to exact
manifold in $\mathbb{R}^4$, but also the approximation of the normal
direction to $\tilde{M}$. In the numerical calculations shown here, we
used a piecewise linear description of the surface of the sphere, and
hence we only expect first order convergence for $\MM{u}$, even if
degree $k>1$ is chosen. As shown in Figure \ref{fig:errs}, we do
indeed obtain first order convergence for $k=0$. For $k=1$, we obtain
second order convergence for $h$ but only first order convergence for
$\MM{u}$. This convergence demonstrates that our methodology produces
convergent solutions under the shallow atmosphere approximation.

\section{Conclusions}

In this paper we introduced a method for implementing the shallow
atmosphere approximation with three dimensional finite element
methods, by making use of a transformation from a three dimensional
manifold embedded in $\mathbb{R}^4$. This can be implemented in a
three dimensional finite element code by defining a transformation to
a discontinuous coordinate field in $\mathbb{R}^3$. This methodology
was demonstrated by numerical convergence tests for a prototype
elliptic problem, implemented using the Firedrake software framework.

\label{sec:summary}

\bibliography{sam-fem}
\end{document}